\documentclass[12pt,twoside]{amsart}

\date{}
\allowdisplaybreaks[4] \footskip=15pt
\renewcommand{\uppercasenonmath}[1]{}

\numberwithin{equation}{section} \theoremstyle{plain}
\newtheorem*{thm*}{Main Theorem}
\newtheorem{thm}{Theorem}[section]
\newtheorem{cor}[thm]{Corollary}
\newtheorem*{cor*}{Corollary}
\newtheorem{lem}[thm]{Lemma}
\newtheorem*{lem*}{Lemma}
\newtheorem{prop}[thm]{Proposition}
\newtheorem*{prop*}{Proposition}
\newtheorem{que}[thm]{Question}
\newtheorem*{que*}{Question}
\newtheorem{rem}[thm]{Remark}
\newtheorem*{rem*}{Remark}
\newtheorem{exa}[thm]{Example}
\newtheorem*{exa*}{Example}
\newtheorem{df}[thm]{Definition}
\newtheorem*{df*}{Definition}

\newtheorem*{conj*}{Conjecture}
\newtheorem*{ack*}{ACKNOWLEDGEMENTS}

\newcommand{\pf}{\noindent\begin {proof}}
\newcommand{\epf}{\end{proof}}

\begin{document}
\title{SMALL INJECTIVE RINGS}
\author{L. Shen and J.L. Chen
}

\address{Department of Mathematics, Southeast University, Nanjing, 210096,
P.R.China }
\address{Department of Mathematics, Southeast University, Nanjing, 210096,
P.R.China }
\date{May 19th, 2005}
 \maketitle \baselineskip=20pt
\begin{abstract}
 Let $R$ be a ring, a right ideal $I$ of $R$ is called small if for
every proper right ideal $K$ of $R$, $I+K\neq R$. A ring $R$ is
called right small injective if every homomorphism from a small
right ideal to $R_{R}$ can be extended to an $R$-homomorphism from
$R_{R}$ to $R_{R}$. Properties of small injective rings are
explored and several new characterizations are given for $QF$
rings and $PF$ rings, respectively.

\end{abstract}

\bigskip

\bigskip
\section{Introduction}
Throughout this paper rings are associative with identity. For a
subset $X$ of a ring $R$, the left annihilator of $X$ in $R$ is
${\bf l}(X)=\{r\in R: rx=0$ for all $x\in X\}$. For any $a\in R$,
we write ${\bf l}(a)$ for ${\bf l}(\{a\})$. Right annihilators are
defined analogously.  We write $J=J(R)$,  $ S_r$ and $ S_l$ for
the Jacobson radical, the right socle and the left socle of $R$,
respectively. $I\subseteq ^{ess}R_{R}$ means $I$ is an essential
right ideal. $f=c\cdot (c\in R)$ means $f$ is a map multiplied by
$c$ on the left side. For a right ideal $I$ of $R$
, we write $I_{n}$ for the set of all $n\times 1$ matrices over $I$.\\
\indent In this article,  the definition of small injective ring
is introduced. Several relations between small injectivity and
other injectivities ( self-injectivity, simple injectivity,
$F$-injectivity and $FP$-injectivity) are given. A main theorem of
Yousif and Zhou \cite[Theorem 2.11]{YZ04} is greatly simplified
and improved by Theorem 3.4. It is well known that a ring $R$ is
quasi-Frobenius (or $QF$) if and only if $R$ is left or right
artinian and left or right self-injective. $R$ is right $PF$ if it
is a semiperfect, right self-injective ring with $S_{r}\subseteq
^{ess}R_{R}$. Under small injective condition, we give some new
characterizations of $QF$ rings and right $PF$ rings,
respectively. The following conditions are proved to be equivalent
in Theorem 3.8: (1) $R$ is $QF$. (2) $R$ is right (or left)
perfect, right and left small injective. (3) $R$ is a semilocal
and right small injective ring with $ACC$ (or $DCC$) on right
annihilators.  (4) $R$ is a right small injective ring with $ACC$
on right annihilators in which $S_{r}\subseteq ^{ess}R_{R}$. (5)
$R$ is a semiregular and right small injective ring with $ACC$ on
right annihilators. It is also proved (see Theorem 3.7) that $R$
is right $PF$ if and only if $R$ is a semilocal, right small
injective ring with $S_{r}\subseteq ^{ess}R_{R}$. Several known
results are generalized as
corollaries.\\

\section{Definitions and examples}
\bigskip
\begin{df}\label{def 2.1}
{\rm A  module $M_{R}$ is called small injective if every
homomorphism from a  small right ideal to $M_{R}$ can be extended
 to an $R$-homomorphism from  $R_{R}$ to $M_{R}$. The left side can be
 defined similarly. A ring $R$ is called right small injective
 if it is  small injective as a right $R$-module.  $R$ is called small injective if it is left
 and right small injective.}
\end{df}
Since $J$ is the sum of all small right (or left ) ideals of a
ring $R$, we have
\begin{exa}\label{exa 2.2}
{\rm Every semiprimitive ring (that is $J$=0) is right and left
small injective.}
\end{exa}

\begin{prop}\label{prop 2.3}
The following are equivalent:\\
{\rm (1)} R is semiprimitive.\\
{\rm (2)} Every right {\rm(}or left{\rm )} R-module is
 small injective.\\
 {\rm (3)} Every principal right {\rm (}or left{\rm )} ideal is
 small injective.\\
\end{prop}
\begin{proof}We only prove the right side, the left side is analogously.  It is obvious that $(1)\Rightarrow (2)\Rightarrow (3)$. Now
we  assume (3), if $J$ is nonzero, then there exists a nonzero
small right ideal $xR$ which  is small injective. It is clear that
the inclusion map from $xR$ to $R_{R}$ is split. Thus $xR$ is a
direct summand of $R$, which is a contradiction.
\end{proof}
\begin{prop}\label{prop 2.4}
A direct product of right R-modules $M=\prod$$ M_{i}$ is small
injective if and only if each $ M_{i}$ is small injective.
\end{prop}
\begin{proof}
By definition.
\end{proof}
\begin{prop}\label{prop 2.5}
A direct product of rings R = $\prod_{i\in I} R_{i}$ is right
small injective if and only if $ R_{i}$ is right small injective,
$\forall i\in I$.
\end{prop}
\begin{proof}
Let $\pi_{i}$ and $\iota_{i}$ be the  $i$th projection and the
$i$th inclusion canonically, $i\in I$. If $R$ is right small
injective, for each $i$, suppose $f_{i} : T_{i}\rightarrow R_{i}$
is $R_{i}$-linear where $T_{i}$ is a small right ideal of $R_{i}$.
Then the map $0\times \cdots\times T_{i}\times \cdots\times
0\rightarrow 0\times \cdots\times R_{i}\times \cdots\times 0$
given by $(0, \cdots,t_{i}, \cdots,0)\longmapsto (0,\cdots, f
_{i}(t_{i}), \cdots, 0)$ is $R$-linear with $0\times \cdots\times
T_{i}\times \cdots\times 0$ a small right ideal of $R$, so it has
the form $c\cdot$ where $c\in R$. Thus $f_{i}=\pi_{i}(c)\cdot$.
Conversely, let $\gamma : T\rightarrow R$ be $R$-linear, where $T$
is a small right ideal of $R$. Write $T_{i}=\{x\in R_{i}~|~
\iota_{i}(x)\in T\}$, it is clear that $T_{i}$ is also a small
right ideal of $R_{i}$, $\forall i\in I$. Now define $\gamma_{i}:
T_{i}\rightarrow R_{i}$ by
$\gamma_{i}(x)=\pi_{i}\gamma(\iota_{i}(x))$, $x\in T_{i}, \forall
i\in I$. Since $R_{i}$ is right small injective,
$\gamma_{i}=c_{i}\cdot$, $\forall i\in I$. Thus for each $\bar
t=\langle t_{i}\rangle \in T$, write $\gamma(\bar t)=\bar
s=\langle s_{i}\rangle$. Since $T$ is a right ideal of $R$,
$t_{i}\in T_{i}, \forall i\in I$. Thus $s_{i}=\pi_{i}(\bar
s\cdot\iota_{i}(1_{i}))=\pi_{i}(\gamma(\bar
t)\cdot\iota_{i}(1_{i}))=\pi_{i}\gamma(\bar
t\cdot\iota_{i}(1_{i}))=\pi_{i}\gamma(\iota_{i}(t_{i}))=\gamma_{i}(t_{i})=c_{i}t_{i}$,
whence $\bar s=\langle c_{i}\rangle\cdot\bar t $. So $R$ is right
small injective.
\end{proof}
\begin{exa}\label{exa 2.6}
    {\rm Every right self-injective ring is right small injective. But
    the converse is not true. For example, the ring of integers
    $\mathbb{Z}$ is a semiprimitive ring but not a self-injective
    ring.}

\end{exa}
\begin{exa}\label{exa 2.7}
{\rm The condition that $R$ is right small injective can not imply
that $R$ is left small injective. In \cite{O84}, Osofsky
constructed a ring $R$ which is semiperfect, right self-injective,
but not left self-injective. Then by Theorem 3.4, such ring is
right small injective but not left small injective.}
\end{exa}
\begin{prop}\label{prop 2.8}
Let R be right small injective. If $e^{2}=e\in R$ satisfies ReR=R,
then eRe is right small injective.
\end{prop}
\begin{proof}
Let $S=eRe$ and $\theta : T\longrightarrow S$ be an $S$-linear
map, where $T$ is a small right ideal of $S$. Define
$\bar{\theta}:TR\longrightarrow R_{R}$ by $\bar{\theta}(\sum
t_{i}r_{i})=\sum \theta (t_{i})r_{i}, t_{i}\in T$. Now we prove
that $\bar{\theta}$ is well defined. Let $\sum t_{i}r_{i}=0$. If
$r\in R$, we get $0=\sum t_{i}r_{i}re=\sum t_{i}(er_{i}re)$,
whence $0=\sum \theta (t_{i})(er_{i}re)=[\sum \theta
(t_{i})r_{i}]re$. Since $ReR=R$, $\sum \theta (t_{i})r_{i}=0$. So
$\bar{\theta}$ is well defined. As $J(eRe)=eJe$,  $TR$ is a small
right ideal of $R$. Hence $\bar{\theta}=c\cdot$, where $c\in R$.
Then $\forall t\in T$,
$\theta(t)=e\theta(t)=e\bar\theta(t)=ect=(ec)et=(ece)t$. It
follows that $\theta=(ece)\cdot$, as required.
\end{proof}
\begin{rem}\label{rem 2.9}
{\rm The condition that $ReR=R$ in the above proposition is
necessary. For example (\cite [Example 9]{K95}), let $R$ be the
algebra of matrices, over a field $k$, of the form
\medskip
\centerline{$R=\left[
\begin{array}{ccccccc}
a & x & 0 & 0& 0 & 0 \\
0 & b & 0 & 0& 0 & 0 \\
0 & 0 & c & y& 0 & 0 \\
0 & 0 & 0 & a& 0 & 0 \\
0 & 0 & 0 & 0& b & z \\
0 & 0 & 0 & 0& 0 & c \\
\end{array}
\right]$.} Set $e=e_{11}+e_{22}+e_{44}+e_{55}$, where $e_{ii}$ are
matrix units. It is clear that $e$ is an idempotent of $R$ such
that $ReR\neq R$. Then $R$ is right small injective, but $eRe$ is
not right small injective.}

\end{rem}
\begin{proof}
By \cite [Example 9]{K95}, $R$ is a $QF$ ring and $eRe$ is not a
$QF$ ring. Since $R$ is $QF$,  $R$ is right small injective, and
$eRe$ is artinian. If $eRe$ is right small injective, then $eRe$
is $QF$ by Theorem 3.8, which is a contradiction.
\end{proof}
\begin{que}\label{que 2.10}
Is right small injectivity  a Morita invariant?
\end{que}
The method in the proof of the following theorem is owing to
\cite[Theorem 1]{NPY00}
\begin{thm}\label{thm 2.11}
The following are equivalent for a ring R and an integer
$n\geq1$:\\
{\rm (1)} M$_{n}${\rm(}R{\rm)} is right small injective.\\
{\rm (2)} For each right R-submodule T of J$_{n}$, every R-linear
map $\gamma$:
T$\rightarrow$ R can be extended to R$_{n}$$\rightarrow$ R.\\
{\rm (3)} For each right R-submodule T of J$_{n}$, every R-linear
map $\gamma$: T$\rightarrow$ R$_{n}$ can be extended to
R$_{n}$$\rightarrow$ R$_{n}$.
\end{thm}
\begin{proof}
We prove for $n$=2, the others are analogous. It is well known that $J(M_{n}(R))=M_{n}(J), \forall n\geq 1$.\\
(1)$\Rightarrow$(2). Given $\gamma$: $T\rightarrow R$ where
$T\subseteq J_{2}$, consider the small right ideal $\bar T$=$[T~
T]=\{[\alpha~\beta]|\alpha, \beta\in T\}$ of $M_{2}(R)$. The map
$\bar\gamma$ : $\bar T$$\rightarrow$$M_{2}(R)$ defined by
$\bar\gamma([\alpha~\beta])$=$\left[
                         \begin{array}{cc}
                            \gamma(\alpha)&\gamma(\beta)\\
                            0&0
                            \end{array}\right]$is $M_{2}(R)$-linear. By (1) we have
$\bar\gamma=C\cdot$ where $C\in M_{2}(R)$, so $\gamma=\bar C\cdot$
where $\bar C$ is the first row of
$C$. Hence $\gamma$ can be extended from $R_{2}$ to $R$.\\
(2)$\Rightarrow$(3). Given (2), consider $\gamma$: $T\rightarrow
R_{2}$ where $T\subseteq J_{2}$. Let $\pi_{i}: R_{2}\rightarrow R$
be the $i$th projection, then (2) provides an $R$-linear map
$\gamma_{i}: R_{2}\rightarrow R$ extending $\pi_{i}\circ\gamma$,
$i=1, 2$. Thus $\bar\gamma:R_{2}\rightarrow R_{2}$ extends
$\gamma$ where $\bar\gamma(\bar x)=[\gamma_{1}(\bar x)~
\gamma_{2}(\bar x)]^{T}$
for all $\bar x\in R_{2}$.\\
(3)$\Rightarrow$(1). Write $S=M_{2}(R)$, consider $\gamma$:
$T\rightarrow S_{S}$ where $T$ is a small right ideal of $S$. Then
$T$=$[T_{0}~ T_{0}]$ where $T_{0}=\{\bar x\in J_{2} ~|~ [\bar x
~0]\in T\}$ is a right $R$-submodule of $J_{2}$. Moreover, the
$S$-linearity shows that $\gamma[\bar x~0]=[\bar y~0]$ for some
$\bar y\in R_{2}$, and writing $\bar y=\gamma_{0}(\bar x)$ yields
an $R$-linear map $\gamma_{0}$: $T_{0}\rightarrow R_{2}$ such that
$\gamma[\bar x~0]=[\gamma_{0}(\bar x)~0]$ for all $\bar x\in
T_{0}$. Then $\gamma_{0}$ extends to an $R$-linear map $\bar
\gamma: R_{2}\rightarrow R_{2}$ by (3). Hence $\gamma_{0}=C\cdot$
for some $C\in S$. If $[\bar x~\bar y]\in T$ it follows that
$\gamma([\bar x~\bar y])=\gamma([\bar x~0]+[\bar
y~0]\left[\begin{array}{cc}
                            0&1\\
                            0&0
                            \end{array}\right])=[\gamma_{0}(\bar x)~0]+[\gamma_{0}(\bar y)~0]\left[
                         \begin{array}{cc}
                            0&1\\
                            0&0
                            \end{array}\right]=[C\bar x~C\bar y]=C[\bar x~\bar y]$, which shows
$\gamma=C\cdot$.
\end{proof}

A ring $R$ is called right mininjective if every $R$-homomorphism
from a minimal right ideal of $R$ to $R_{R}$ can be extended from
$R_{R}$ to $R_{R}$. $R$ is called left minannihilator if every
minimal left ideal is a left annihilator.
\begin{prop}\label{prop 2.12}
If R is right small injective, then \\
{\rm(1)} R is right mininjective.\\
{\rm(2)} ${\bf l}(I\cap L)={\bf l}(I)+{\bf l}(L)$, for any small
right ideals I and L of R.\\
{\rm(3)} Every principal small left  ideal of R is a left
annihilator, so R is
left minannihilator. \\
 {\rm(4)} ${\bf l}(bR\cap {\bf
r}(a))={\bf l}(b)+Ra$, $\forall b\in R, a\in J$.

\end{prop}
\begin{proof}
For (1), since every minimal one-sided ideal of $R$ is either
nilpotent or a direct summand of $R$ (see \cite[(10.22) Brauer's
Lemma]{L91}), each right small injective ring is right
mininjective. But the converse is not true ( see Example 2.13).\\
(2) and (3) by the similar proof of \cite[Lemma 30.9 ]{AF92}.\\
(4) see \cite[Lemma 1.2]{YZ04}.
\end{proof}
\begin{exa}\label {exa 2.13}
{\rm({\rm The Bj\"{o}rk Example} \cite [Example 2.5]{NY03})} {\rm
A right mininjective ring may not be right small injective. Let
$F$ be a field and assume that $a$ $\mapsto \bar{a}$ is an
isomorphism $F$ $\mapsto \bar{F}\subseteq F$, where the subfield
$\bar{F}\neq F$. Let $R$ denote the left vector space on basis
\{1,$t$\}, and make $R$ into an $F$-algebra by defining $t^{2}$=0
and $ta$=$\bar{a}$t for all $a\in F$. Then $R$ is a right
mininjective ring but not a right small injective ring.}
\end{exa}

\begin{proof}
It is mentioned in  \cite [Example 2.5]{NY03} that $R$  is
semiprimary, right mininjective but not left mininjective. If $R$
is right small injective, then it is right self-injective by
Theorem 3.4. Thus it is a right $PF$ ring, which shows $R$ is left
mininjective by \cite[Theorem 5.57]{NY03}, a contradiction.
\end{proof}
A ring $R$ is called right Kasch if every simple right $R$-module
can embed into $R_{R}$.
\begin{prop}\label {prop 2.14}
If R is right small injective and right Kasch, then\\
{\rm(1)} ${\bf r}{\bf l}(I)=I$ for every small right ideal I of R.\\
{\rm(2)} The map $\theta :T\rightarrow {\bf l}(T)$ from the set of
maximal right ideals T of R to the set of minimal left ideals of R
is a bijection. And the inverse map is given by K$\longmapsto {\bf
r}(K)$, where K is a minimal left ideal of R.\\
 {\rm(3)} For $k\in R$, Rk is minimal if and only if kR is minimal,
in particular $S_{r}=S_{l}$.\\

\end{prop}
\begin{proof}
(1) is by \cite[Lemma 2.4 (3)]{YZ04}.\\
(2) is informed by \cite[Theorem 2.32 (b)]{NY03} and Proposition
2.12 (3).\\
For (3), if $Rk$ is minimal, then ${\bf r}(k)$ is maximal by (2)
which shows $kR$ is also minimal. Conversely, if $kR$ is minimal,
then $Rk$ is minimal by \cite[Theorem 2.21 (a)]{NY03}.
\end{proof}

\section{some relations between small injectivity and other
injectivities}
\begin{lem}\label{lem 3.1}
Let $R$ be a semilocal ring and I  a right ideal of R, then every
homomorphism from a right ideal to $I$ can be extended to an
endomorphism of $R_{R}$ if and only if every homomorphism from a
small right ideal to I can be extended to an endomorphism of
$R_{R}$ .
\end{lem}
\begin{proof}
             {\rm(i)} $``\Longrightarrow"$ is obvious.
\\\indent\indent$~$ {\rm(ii)}$``\Longleftarrow"$ Let $f$ be a
homomorphism from a right ideal $K$ of $R$ to $I$. Since $R$ is
semilocal, there exists a right ideal $L$ of $R$ such that $K+L=R$
and $K\cap L\subseteq J$(see \cite[Corollary 3.2]{L99}). Thus
$K\cap L$ is small and there exists an endomorphism $g$ of $R_{R}$
such that $g(x)=f(x)$, $\forall x\in K\cap L$. Define $F$:
$R_{R}\longrightarrow R_{R}$ such that for any $x=k+l, k\in K,
l\in L$, $F(x)=f(k)+g(l)$. Now we prove that $F$ is well defined.
If $k_{1}+l_{1}=k_{2}+l_{2}, k_{i}\in K, l_{i}\in L, i=1,2$, then
$k_{1}-k_{2}=l_{2}-l_{1}\in K\cap L$. Hence
$f(k_{1}-k_{2})=g(l_{2}-l_{1})$, which shows
$F(k_{1}+l_{1})=F(k_{2}+l_{2})$. Thus $F$ is an endomorphism of
$R_{R}$ such that $F_{\mid K}$=$f$.
\end{proof}
A right ideal $I$ of $R$ is said to lie over a summand of $R_{R}$
if there exists a direct decomposition $R_{R}=P_{R}\oplus Q_{R}$
with $P\subseteq I$ and $Q\cap I$ is small in $R$. In this case,
$I=P\oplus (Q\cap I)$.
\begin{lem}\label{lem 3.2}
Let $R$ be a  ring and I  a right ideal of R. If every
{\rm(}m-generated{\rm)} right ideal lies over a summand of
R$_{R}$, then every homomorphism from an {\rm(}m-generated{\rm)}
small right ideal to I can be extended to an endomorphism of R if
and only if every homomorphism from an {\rm(}m-generated{\rm)}
right ideal to I can be extended to an endomorphism of R.

\end{lem}
\begin{proof}
Let $K$ be any ($m$-generated) right ideal of $R$ and $f$ a
homomorphism from $K$ to $I$. Since $K$ lies over a summand of
$R_{R}$, there exists an idempotent $e^{2}=e\in R$ such that
$K=eR\oplus L$, where $L\subseteq J$ is  an ($m$-generated) small
right ideal. Now we show that $K=eR\oplus (1-e)L$. It is obvious
that $eR+(1-e)L$ is a direct sum. If $x\in eR\oplus L$, then
$x=a+b$, for some $a\in eR$, $b\in L$. Thus $x=a+eb+(1-e)b\in
eR\oplus (1-e)L$. On the other hand, if $y\in eR\oplus (1-e)L$,
write $y=c+(1-e)d$ where $c\in eR$, $d\in L$. Hence
$y=c+(1-e)d=(c-ed)+d\in eR\oplus L$. As $J$ is an ideal which is a
small right ideal, so $(1-e)L\subseteq J$ is also an
($m$-generated) small right ideal. Since every homomorphism from
an ($m$-generated) small right ideal to $I$ can be extended to an
endomorphism of $R$, there exists a homomorphism $g$ from $R_{R}$
to $R_{R}$ such that $g_{\mid(1-e)L}$=$f_{\mid(1-e)L}$. Now we
define $F$ from $R_{R}$ to $R_{R}$ such that
$F(x)=f(ex)+g((1-e)x)$, whence $F$ is a well defined homomorphism.
Then for every $x=a+b\in K=eR\oplus (1-e)L$ where $a\in eR,~ b\in
(1-e)L$, $F(x)=f(ex)+g((1-e)x)=f(a)+g(b)=f(a)+f(b)=f(a+b)=f(x)$,
which
shows that $F_{\mid K}$=$f$. \\
\indent The converse is obvious.

\end{proof}

 Let $I$, $K$
be two right ideals of $R$ and $m\geq1$. $R$ is called right
($I,K$)-$m$-injective (see \cite{YZ04}) if, for any $m$-generated
right ideal $U\subseteq I$ and any $R$-homomorphism from $U_{R}$
to $K_{R}$ can be extended from $R_{R}$ to $R_{R}$. A ring $R$ is
called right simple injective (right simple $J$-injective) if for
any homomorphism from a right ideal (small right ideal) of $R$ to
$R_{R}$ with simple image can be extended from $R_{R}$ to $R_{R}$.
$R$ is called right ($m$-injective) $F$-injective if every
homomorphism from an ($m$-generated) finitely generated right
ideal of $R$ to $R_{R}$ can be extended from $R_{R}$ to $R_{R}$.
$R$ is called right $(I,K)$-$FP$-injective if, for any $n\geq 1$
and any finitely generated $R$-submodule $N$ of $I_{n}$, every
$R$-homomorphism $f: N\rightarrow K$ can be extended to an
$R$-homomorphism $g: R_{n}\rightarrow R$. $R$ is right
$FP$-injective if $R$ is right $(R,R)$-$FP$-injective.
\begin{lem}\label{lem 3.3}
\cite[Lemma 1.3]{YZ04}The following are equivalent for
two right ideals I and K of R:\\
{\rm(1)} R is right {\rm(}I,K{\rm)}-FP-injective.\\
{\rm(2)} M$_{n}${\rm(}R{\rm)} is right
{\rm(}M$_{n}{\rm(}I{\rm)},M_{n}{\rm(}K{\rm)}${\rm)}-1-injective
for every n $\geq$ 1.
\end{lem}
A ring $R$ is called semiregular if $R/J$ is von Neuman regular
and idempotents lift modulo $J$.
 \begin{thm} \label{thm 3.4}
Let R be a ring, we have \\
{\rm(1)} If R is semilocal, then R is right self-injective if and
only if R is right
small injective.\\
{\rm(2)} If R is semilocal, then R is right simple injective if
and only
if R is right simple J-injective.\\
 {\rm(3)} If R is semiregular, I is a right ideal of R. Then R is right {\rm(}J,I{\rm)}-m-injective if and only
if R is right {\rm(}R,I{\rm)}-m-injective. In particular, R is
right {\rm(}J,S$_{r}${\rm)}-m-injective if and only if R is right
{\rm(}R,S$_{r}${\rm)}-m-injective, R is right
{\rm(}J,R{\rm)}-m-injective if and only if R is right m-injective,
R is right F-injective if and only if R is right
{\rm(}J,R{\rm)}-k-injective, $\forall k\geq1$.
\\
{\rm(4)} If R is semiregular, then R is right
{\rm(}J,R{\rm)}-FP-injective if and only if R is right
FP-injective.
\end{thm}
\begin{proof}
(1) and (2) by Lemma 3.1. For (3), since $R$ is semiregular, every
finitely generated right (or left) ideal of $R$ lies over a
summand of $R$ (see \cite [Theorem 2.9]{N76}). Then (3) is clear
by Lemma 3.2.
 Since semiregularity is a Morita invariant (see \cite[Corollary
2.8]{N76}) and  $J(M_{n}(R))$=$M_{n}(J)$, (4) is followed by (3)
and Lemma 3.3 .

\end{proof}
Since a semiperfect ring is   both semilocal and semiregular, we
have

 \begin{cor} \label{cor 3.5}
Let R be a semiperfect ring, we have \\
{\rm (1)} R is right self-injective if and only if R is right
small injective.\\
{\rm (2)} R is right simple injective if and only
if R is right simple J-injective.\\
 {\rm (3)} Let I be a right ideal of R. Then R is right {\rm(}J,I{\rm)}-m-injective if and only
if R is right {\rm(}R,I{\rm)}-m-injective. In particular, R is
right {\rm(}J,S$_{r}${\rm)}-m-injective if and only if R is right
{\rm(}R,S$_{r}${\rm)}-m-injective, R is right
{\rm(}J,R{\rm)}-m-injective if and only if R is right m-injective,
R is right F-injective if and only if R is right
{\rm(}J,R{\rm)}-k-injective, $\forall k\geq1$.
\\
{\rm (4)} R is right {\rm(}J,R{\rm)}-FP-injective if and only if R
is right FP-injective.
\end{cor}

 \begin{cor}\label {cor 3.6}
 \cite[Theorem 2.11]{YZ04} Let R be a
semiperfect ring with $S_{r}\subseteq^{ess}R_{R}$, then\\
{\rm(1)} If R is right
{\rm(}J,S$_{r}${\rm)}-{\rm(}m+1{\rm)}-injective, then R is right
{\rm(}R,S$_{r}${\rm)}-m-injective.\\
 {\rm(2)} If R is right {\rm(}J,R{\rm)}-{\rm(}m+1{\rm)}-injective, then R is right
m-injective.\\
{\rm(3)} If R is right simple J-injective, then R is right simple
injective.\\
{\rm(4)} If every  homomorphism from a small right ideal of R to R
can be extended to an R-homomorphism from $R_{R}$ to $R_{R}$, then
R is right self-injective.

\end{cor}

\begin{thm}\label{thm 3.7}
The following are equivalent:\\
{\rm(1)} R is right PF.\\
{\rm(2)} $R$ is a semilocal, right small injective ring with
$S_{r}\subseteq ^{ess}R_{R}$.

\end{thm}
\begin{proof}
It is obvious that (1) implies (2). By Theorem 3.4, the ring
satisfying (2) is right self-injective, which shows that
idempotents can be lifted modulo $J$. Hence $R$ is semiperfect.
Thus (2) implies (1).
\end{proof}
\begin{thm}\label{thm 3.8}
The following are equivalent:\\
{\rm(1)} R is QF; \\
{\rm(2)} R is right (or left) perfect, right and left small injective; \\
{\rm (3)} R is a semilocal and right small injective ring with
$ACC$ (or $DCC$) on right annihilators. \\
{\rm(4)} R is a right small injective ring with ACC on right
annihilators in which $S_{r}\subseteq ^{ess}R_{R}$.
\\
{\rm (5)} R is a semiregular and right small injective ring with
$ACC$  on right annihilators. \\
\end{thm}
\begin{proof}
It is obvious that (1)$\Rightarrow$ (2), (3), (4) and (5). It is
proved in \cite[Lemma 2.11]{SC} that a right mininjective ring
with $ACC$ on right annihilators in which $S_{r}\subseteq
^{ess}R_{R}$ is semiprimary. Then the rings in (2)-(4) are all
semilocal rings. Hence they are all right self-injective by
Theorem 3.4. Thus (2), (3) and (4) are clear by \cite[Theorem 2.3,
Theorem 4.1(b)]{FH02}. For (5), $R$ is right $F$-injective by
Theorem 3.4. Then (5) is implied by \cite[Corollary 3]{R75}.
\end{proof}
\section{extensions of small injective rings}
Given a ring $R$ and a bimodule $_{R}V_{R}$, the trivial extension
of $R$ by $V$ is the ring $S=R\propto V=\{(r,v): r\in R, v\in V\}$
with the usual addition and multiplication
$(r,v)(r^{'},v^{'})=(rr^{'},rv^{'}+vr^{'})$. In fact, $R\propto V$
is isomorphic to the ring of all matrices $\left[
                         \begin{array}{cc}
                            r&v\\
                            0&r
                            \end{array}\right]$
                                                     where $r\in R$ and
$v\in V$ with the usual  matrix operations. For convenience, we
let $I\propto V=\{(r,v): r\in I, v\in V\}$ where $I$ is a subset
of $R$. Clearly, $V$ is an ideal of $S$, $V^{2}=0$ and $S/V\cong
R$. The Jacobson radical of $S$ is $J(R)\propto V$.
\begin{prop}\label {prop 4.1}
Let $S=R\propto V$, where R is a ring and $_{R}V_{R}$ a bimodule
of R. If S is right small injective, then V is self-injective as a
right R-module and R = End $V_{R}$ canonically.
\end{prop}
\begin{proof}
Let $K_{R}$ be a right $R$-submodule of $V$ and $f$ a right
$R$-homomorphism from $K$ to $V$. Then $\{\left[
                         \begin{array}{cc}
                            0&a\\
                            0&0
                            \end{array}\right]:a\in K\}$ is a small
right ideal of $S$. Now we \\ define $g$ from $\{\left[
                         \begin{array}{cc}
                            0&a\\
                            0&0
                            \end{array}\right]:a\in K\}$  to $S$ such that
$g(
                     \left[\begin{array}{cc}
                            0&a\\
                            0&0
                            \end{array}\right])$=$\left[
                         \begin{array}{cc}
                            0&f(a)\\
                            0&0
                            \end{array}\right],
a\in K$. It is clear that $g$ is  right $S$-linear. Since $S$ is
right small injective, there \\ exists $\left[
                         \begin{array}{cc}
                            b&c\\
                            0&b
                            \end{array}\right]\in S$ satisfying $g(
                     \left[\begin{array}{cc}
                            0&a\\
                            0&0
                            \end{array}\right])$=$\left[
                         \begin{array}{cc}
                            0&f(a)\\
                            0&0
                            \end{array}\right]
= \left[\begin{array}{cc}
                            b&c\\
                            0&b
                            \end{array}\right]
                            \left[\begin{array}{cc}
                            0&a\\
                            0&0
                            \end{array}\right]$=\\$\left[\begin{array}{cc}
                            0&ba\\
                            0&0
                            \end{array}\right]$. Thus $f(a)=ba, \forall a\in K$, which shows that $V_{R}$ is
right self-injective and  $R$ = End $V_{R}$ canonically.
\end{proof}
\begin{cor}\label{cor 4.2}
Let $S=R\propto V$, where R is a ring and $_{R}V_{R}$ a bimodule
of R. If V is injective as a right R-module. Then the following
are
equivalent:\\
{\rm(1)} S is right self-injective;\\
{\rm(2)} S is right small injective;\\
{\rm(3)}  R = End $V_{R}$ canonically;
\end{cor}
\begin{proof}
 By \cite[Theorem 2]{F79} and the above proposition.
\end{proof}
\begin{cor}\label{cor 4.3}
Let R be a ring and  $S=R\propto R$. Then the following are
equivalent:\\
{\rm(1)} R is right self-injective;\\
{\rm(2)} S is right self-injective;\\
{\rm(3)} S is right small injective;\\

\end{cor}

\begin{rem}\label{rem 4.4}
{\rm In the above corollary, $R$ is right small injective can not
implies that $S$ is right small injective. For example, let
$S$=$\mathbb{Z}$ $\propto$ $\mathbb{Z}$.  Since $\mathbb{Z}$ is
semiprimitive, $\mathbb{Z}$ is right small injective. But $S$ can
not be right small injective. If not, by the above corollary, we
have that $\mathbb{Z}$ is self-injective, which is a
contradiction.}
\end{rem}
If $R$ and $S$ are rings and $_{R}V_{S}$ is a bimodule, the formal
triangular matrix ring of $R$ and $S$ by $V$ is the ring $U=\left[
                         \begin{array}{cc}
                            R&V\\
                            0&S
                            \end{array}\right]$
 of all matrices$\left[
                         \begin{array}{cc}
                            r&v\\
                            0&s
                            \end{array}\right]$ where $r\in R, s\in S, v\in V$ with the usual matrix
                            operations. It is clear that
                            $J(U)=$$\left[
                         \begin{array}{cc}
                            J(R)&V\\
                            0&J(S)
                            \end{array}\right]$.
 \begin{prop}\label{prop 4.5}
 Assume that $U=\left[
                         \begin{array}{cc}
                            R&V\\
                            0&S
                            \end{array}\right]$ is right small
                            injective. Then:\\
 {\rm(1)} S is right small injective.\\
 {\rm(2)} For every right S-submodule $K$ of $V_{S}$, {\rm Hom}$_{S}(K,
 S)=0$.\\
 {\rm(3)} If $\gamma : K_{S}\rightarrow V_{S}$ is S-linear, where K is an S-submodule of
 $ V_{S}$, then $\gamma =r\cdot$ for some $r\in R$.\\
 {\rm(4)} $X_{S}\subseteq\left[\begin{array}{c}
                            V_{S}\\
                            J(S)
                            \end{array}\right]$ and $\theta: X_{S}\rightarrow
 \left[\begin{array}{c}
                            V_{S}\\
                            S
                            \end{array}\right]$ is S-linear, then $\theta=C\cdot$
 for some $C\in U$.

 \end{prop}
 \begin{proof}
 (1). If $\gamma : T\rightarrow S_{S}$, where $T$ is a small right
 ideal of $S$ and $\gamma$ is $S$-linear. Then $\bar T=\left[\begin{array}{cc}
                            0&0\\
                            0&T
                            \end{array}\right]$ is a small right
                            ideal of $U$ and $\bar\gamma: \bar
                            T\rightarrow U$ is $S$-linear by
                            defining $\bar\gamma(\left[\begin{array}{cc}
                            0&0\\
                            0&t
                            \end{array}\right])=\left[\begin{array}{cc}
                            0&0\\
                            0&\gamma(t)
                            \end{array}\right]$. Hence $\bar\gamma=\left[\begin{array}{cc}
                            r&v\\
                            0&s
                            \end{array}\right]\cdot$, for some $r\in R, s\in S, v\in V$.
                            Thus
                            $\gamma=s\cdot$.\\
    (2). Let $\gamma: K_{S}\rightarrow S_{S}$ where $K_{S}$ is a
    right $S$-submodule of $V_{S}$. Then $\bar K=\left[\begin{array}{cc}
                            0&K\\
                            0&0
                            \end{array}\right]$ is a small right
                            ideal of $U$ and $\bar\gamma : \bar
                            K\rightarrow U$ is $U$-linear by
                            defining $\bar\gamma(\left[\begin{array}{cc}
                            0&k\\
                            0&0
                            \end{array}\right])=\left[\begin{array}{cc}
                            0&0\\
                            0&\gamma(k)
                            \end{array}\right]$. Hence $\bar\gamma=\left[\begin{array}{cc}
                            r&v\\
                            0&s
                            \end{array}\right]\cdot$, for some $r\in R, s\in S, v\in V$. So $\left[\begin{array}{cc}
                            0&0\\
                            0&\gamma(k)
                            \end{array}\right]=\left[\begin{array}{cc}
                            r&v\\
                            0&s
                            \end{array}\right]\left[\begin{array}{cc}
                            0&k\\
                            0&0
                            \end{array}\right]=\left[\begin{array}{cc}
                            0&rk\\
                            0&0
                            \end{array}\right]$, which implies
                            $\gamma=0$.\\
 (3) If $\gamma$ is as given, then $\bar\gamma :\bar K\rightarrow U$ is $U$-linear if we define $\bar\gamma(\left[\begin{array}{cc}
                            0&k\\
                            0&0
                            \end{array}\right])=\left[\begin{array}{cc}
                            0&\gamma(k)\\
                            0&0
                            \end{array}\right]$. Hence
                            $\bar\gamma=\left[\begin{array}{cc}
                            r&v\\
                            0&s
                            \end{array}\right]\cdot$, for some $r\in R, s\in S, v\in V$. Hence
                            $\gamma=r\cdot$.\\
 (4). Write $[0~ X]=\{\left[\begin{array}{cc}
                            0&v\\
                            0&s
                            \end{array}\right]|\left[\begin{array}{cc}
                            v\\
                            s
                            \end{array}\right]\in X\}$. Then $[0~
                            X]$ is a small right ideal of $U$ and
                            $\bar\theta :[0~X]\rightarrow U$ is
                            $U$-linear by defining $\bar\theta\left[\begin{array}{cc}
                            0&v\\
                            0&s
                            \end{array}\right]=\left[\begin{array}{cc}
                            0\\
                            0
                            \end{array}\theta(\left[\begin{array}{cc}
                            v\\
                            s
                            \end{array}\right])\right]$. Thus
                            $\bar\theta=C\cdot$ for some $C\in U$, and so is
                            $\theta=C\cdot$.
 \end{proof}
 Let $S=T_{n}(R)$ be the upper triangular matrix ring over a ring
 $R$, $n\geq 2$. We have
 \begin{prop}\label{prop4.6}
 For each $n\geq 2$, $S=T_{n}(R)$ can not be right small
 injective.
 \end{prop}
 \begin{proof}
 If not, we assume $S$ is right small injective. Let $e_{ij}$ be usual matrix units, 1$\leq i, j\leq n$. If 0$\neq
 x\in R$, then $K=e_{1n}xR$ is a small right ideal of $S$. Now we
 define $\gamma : K\rightarrow S$ by $\gamma(e_{1n}xt)=e_{nn}xt, \forall t\in R$. It
 is clear that $\gamma$ is right $S$-linear. Thus $\gamma=C\cdot,$
 for some $C\in S$. Hence $e_{nn}x=Ce_{1n}x$, which shows $x=0$, a
 contradiction.
 \end{proof}
\medskip

\begin{center}
{ACKNOWLEDGEMENTS}
 \end{center}
\indent\indent The research is supported by the National Natural
Science Foundation of China (No.10171011) and the Teaching and
Research Award Program for Outstanding Young Teachers in Higher
Education Institutes of MOE, P.R.C.

\end{document}